\documentclass[leqno,11pt]{article}
\usepackage{amssymb,amsmath, amsthm,latexsym, verbatim, amscd}
\usepackage{graphicx}

\usepackage{color}
\usepackage[hypertex]{hyperref}%

\newtheorem{theorem}{Theorem}[section]
\newtheorem{prop}[theorem]{Proposition}
\newtheorem{lemma}[theorem]{Lemma}
\newtheorem{cor}[theorem]{Corollary}
\newtheorem{definition}[theorem]{Definition}
\newtheorem{example}[theorem]{Example}
\newtheorem{remark}[theorem]{Remark}

\newcommand{\btheorem}{\begin{theorem}}
\newcommand{\etheorem}{\end{theorem}}
\newcommand{\bprop}{\begin{prop}}
\newcommand{\eprop}{\end{prop}}
\newcommand{\blemma}{\begin{lemma}}
\newcommand{\elemma}{\end{lemma}}
\newcommand{\bcor}{\begin{cor}}
\newcommand{\ecor}{\end{cor}}

\newcommand{\bg}{{\mathfrak g }}

\newcommand{\bq}{{\mathfrak q }}

\newcommand{\bC}{{\mathbb C}}
\newcommand{\bZ}{{\mathbb Z}}

\newcommand{\Exterior}{\mathchoice{{\textstyle\bigwedge}}
    {{\bigwedge}}
    {{\textstyle\wedge}}
    {{\scriptstyle\wedge}}}

\begin{document}
\title{\bf Symmetry breaking  
for orthogonal groups and a conjecture by B.~Gross and D.~Prasad}
\author{
Toshiyuki Kobayashi
\vspace{-10mm}
\footnote
{
Kavli IPMU and Graduate School of Mathematical Sciences,
 the University of Tokyo, 3-8-1 Komaba, Tokyo, 153-8914 Japan
\newline
\textit{E-mail address}:
\texttt{toshi@ms.u-tokyo.ac.jp}
}
\,\,
 and 
 Birgit Speh
\vspace{-10mm}
\footnote
{Department of Mathematics, Cornell University, Ithaca, NY 14853-4201, USA
\newline
\textit{E-mail address}:
\texttt{bes12@cornell.edu}}
}

\maketitle

\pagestyle{myheadings}
   \markboth{Symmetry breaking and the Gross--Prasad conjecture}
            {Symmetry breaking and the Gross--Prasad conjecture}

\begin{abstract}
We consider {\bf irreducible unitary representations $A_i$  } of $G=SO(n+1,1)$ with the same infinitesimal character  as the trivial representation  and representations $B_j$ of $H=SO(n,1)$ with the same properties and discuss $H$-equivariant homomorphisms $\operatorname{Hom}_H(A_i,B_j)$.   For tempered representations our results confirm the predictions of conjectures by B.~Gross and D.~Prasad.
\end{abstract}

\section{Introduction}
\label{sec:intro}
A representation $\Pi$ of a group $G$
 defines a representation
 of a subgroup $G'$
 by restriction.  
In general irreducibility 
 is not preserved by the restriction.  
If $G$ is compact
 then the restriction $\Pi|_{G'}$ is isomorphic
 to a direct sum of irreducible finite-dimensional representations $\pi$ of $G'$ with multiplicities  $m(\Pi,\pi)$. 
These multiplicities are studied by using combinatorial techniques.  
We are interested in the case
 where $G$ and $G'$ are (noncompact) real reductive Lie groups.  
Then most irreducible representations $\Pi$ of $G$ are
 infinite-dimensional, 
 and generically the restriction $\Pi|_{G'}$
 is not a direct sum of irreducible representations
 \cite{KInvent98}.  
So we have to consider another notion of multiplicity.

For a continuous representation $\Pi$ of $G$
 on a complete, 
 locally convex topological vector space ${\mathcal{H}}$,
 the space ${\mathcal{H}}^\infty $ of $C^\infty$-vectors of ${\mathcal{H}}$
 is naturally endowed with a Fr{\'e}chet topology,
 and $(\Pi,{\mathcal{H}})$ induces a continuous representation $\Pi^{\infty}$ of $G$
 on ${\mathcal{H}}^\infty$.  
If $\Pi$ is an admissible representation
 of finite length on a Banach space ${\mathcal{H}}$, 
 then the Fr{\'e}chet representation 
 $(\Pi^{\infty}, {\mathcal{H}}^{\infty})$
 depends only 
 on the underlying $({\mathfrak {g}}, K)$-module
 ${\mathcal{H}}_K$, 
 sometimes referred to as an admissible representation
 of moderate growth \cite[Chap.~11]{W}.  
We shall work with these representations
 and write simply $\Pi$ for $\Pi^{\infty}$.  
Given another continuous representation
 $\pi$ of moderate growth
 of a reductive subgroup $G'$, 
 we consider the space of continuous $G'$-intertwining operators 
 ({\it{symmetry breaking operators}})
\[ \operatorname{Hom}_{G'} ({\Pi}|_{G'}, {\pi}) .\] 
The dimension $m(\Pi,\pi) $ of this space 
 yields important information of the restriction of $\Pi $ to $G'$
 and is called the {\it{multiplicity}} of $\pi$
 occurring in the restriction $\Pi|_{G'}$.  
In general,
 $m(\Pi,\pi)$ may be infinite.  
  The criterion in \cite{KO} asserts that the
 multiplicity $m(\Pi,\pi)$ is finite
 for all irreducible representations $\Pi$ of $G$
 and all irreducible representations $\pi$ of $G'$
 if and only if a minimal parabolic subgroup $P'$
 of $G'$ has an open orbit 
 on the real flag variety $G/P$, 
 and that the multiplicity is uniformly bounded
 with respect to $\Pi$ and $\pi$
 if and only if a Borel subgroup of $G_{\mathbb{C}}'$
 has an open orbit
 on the complex flag variety of $G_{\mathbb{C}}$.  

\medskip

We consider in this article the case 
\begin{equation}
\label{eqn:GG}
     (G,G')=(SO(n+1,1), SO(n,1)), 
\end{equation}
 and discuss  symmetry breaking 
 between irreducible unitary representations
 of the groups $G$ and $G'$
 with the same infinitesimal character
 $\rho$
 as the trivial one-dimensional representations.

\medskip
We state our results first in Langlands parameters
 by identifying the representations with the Langlands subquotients of principal series representations induced from finite-dimensional representations of a maximal parabolic subgroup. 
Since these representations also have nontrivial $(\bg,{\mathfrak{k}})$-cohomology we can parametrize them by characters of the Levi of a $\theta$-stable parabolic subalgebras and we proceed
 to state the results in this language. Then we describe the representations as members of Vogan packets and restate the results in this language. 
In the last section we relate our results to the Gross--Prasad conjectures for tempered representations.
 
 \medskip
 Detailed proofs of the results will be published elsewhere.

\section{\bf Classification of symmetry breaking operators}
\label{sec:notations}
\medskip
The main result of this section
 is a classification
 of symmetry breaking operators
 for principal series representations
 induced from exterior tensor representations
 for the pair $(G,G')=(SO(n+1,1),SO(n,1))$.  
Theorem \ref{thm:1.1} extends the scalar case \cite{KS}
 and the case of differential operators \cite{KKP}, 
 and will be used in Section \ref{sec:mainresults}.  

\subsection{Notation for $SO(n+1,1)$}
\label{subsec:II.1}
We first recall the  notation from the Memoir article
\cite{KO}.  

\medskip

Consider the quadratic form 
\begin{equation}
\label{eqn:quad}
      x_0^2 + x_1^2 +\dots +x_{n}^2-x_{n+1}^2
\end{equation}
 of signature $(n+1,1)$. 
We define $G$
 to be the indefinite special orthogonal group $SO(n+1,1)$
 that preserves the quadratic form 
 \eqref{eqn:quad} and the orientation.  
Let $G'$ be the stabilizer
 of the vector $e_{n}={}^t\! (0,0,\cdots,0,1,0)$.  
Then $G' \simeq SO(n,1)$.  
We set 
\begin{alignat}{3}
\label{eqn:K}
K &:=O(n+2) \cap G
&&=\{
\begin{pmatrix}
A &  
\\
  & \det A
\end{pmatrix}
:
A \in O(n+1) \} 
&& \simeq O(n+1), 
\\
K' &:= K \cap G' 
&&=\{
\begin{pmatrix}
B & & 
\\
  & 1 &
\\
  & & \det B
\end{pmatrix}
:
B \in O(n)
\}
&&\simeq O(n).  
\notag
\end{alignat}
Then $K$ and $K'$ are maximal compact subgroups
 of $G$ and $G'$, 
respectively.

Let ${\mathfrak {g}}={\mathfrak {so}}(n+1,1)$
 and ${\mathfrak {g}}'={\mathfrak {so}}(n,1)$
 be the Lie algebras of $G$
 and $G'$, 
 respectively.  
We take a hyperbolic element $H$
 as 
\begin{equation}
H :=
  E_{0,n+1} + E_{n+1,0} \in \mathfrak{g}', 
\end{equation}
and set
\begin{equation*}
{\mathfrak {a}}:={\mathbb{R}}H
\qquad
\text{ and }
A:=\exp {\mathfrak {a}}.  
\end{equation*}
Then the centralizers of $H$ in $G$ and $G'$
 are given by $M A$ and $M' A$, 
 respectively,
 where 
\begin{align*}
M :={}
  & \left\{
    \begin{pmatrix} 
        \varepsilon \\ & A \\ & & \varepsilon
    \end{pmatrix} :
    A \in SO(n), \   \varepsilon = \pm 1
    \right\} 
   & \simeq & SO(n) \times O(1), 
\\
M' :={}
  &  \left\{
    \begin{pmatrix} 
        \varepsilon \\ & B \\ & & 1 \\ & & & \varepsilon
    \end{pmatrix} :
    B \in SO(n-1):   \varepsilon = \pm 1
    \right\}
 &  \simeq & SO(n-1) \times O(1).  
\label{eqn:Mprime}
\end{align*}
We observe that 
 $\operatorname{ad}(H) \in \operatorname{End}_{\mathbb{R}}(\mathfrak{g})$
 has eigenvalues 
 $-1$, $0$, and $+1$.  
Let 
\[
{\mathfrak {g}} ={\mathfrak {n}}_- + ({\mathfrak {m}}+{\mathfrak {a}})+{\mathfrak {n}}_+
\]
be the corresponding eigenspace decomposition, 
 and $P$ a minimal parabolic subgroup  with $P=MAN_+$ its Langlands decomposition. 
We remark
 that $P$ is also a maximal parabolic subgroup of $G$.  
Likewise, 
 $P'=M'AN_+'$ is a compatible Langlands decomposition
 of a minimal (also maximal) parabolic subgroup $P'$ of $G'$ given by 
 \begin{equation}
\label{eqn:Lang}
{\mathfrak {p}}'={\mathfrak {m}}'+{\mathfrak {a}}+{\mathfrak {n}}_+'
 = ({\mathfrak {m}} \cap {\mathfrak {g}}')
  +({\mathfrak {a}} \cap {\mathfrak {g}}')
  +({\mathfrak {n}}_++{\mathfrak {g}}').  
\end{equation}

We note 
 that we have chosen $H \in {\mathfrak {g}}'$
 so that $P'=P \cap G'$
 and we can take a common maximally split abelian subgroup
 $A$ in $P'$ and $P$.  

\medskip
\subsection{Principal series representations of $SO(n+1,1)$}
\label{subsec:II.2}

The character group of $O(1)$
 consists of two characters.  
We write $+$ for the trivial character,
 and $-$ for the nontrivial character.  
Since $M \simeq SO(n) \times O(1)$, 
 any irreducible representation
 of $M$ is the outer tensor product
 of a representation $(\sigma, V)$ of $SO(n)$ and a character 
 $\delta$ of $O(1)$.

Given $(\sigma,V) \in \widehat {SO(n)}$, 
 $\delta \in \{ \pm \} \simeq \widehat{O(1)}$, 
 and a character $e_{\lambda}(\exp (t H))=e^{\lambda t}$
 of $A$ for $\lambda \in {\mathbb{C}}$, 
 we define the (unnormalized) principal series representation
\[
   I_{\delta}(V, \lambda)= \operatorname{Ind}_{P}^{G}
  (V \otimes \delta, \lambda) 
\]
of $G=SO(n+1,1)$ on the Fr{\'e}chet space 
 of smooth maps
 $f \colon G \to V$ 
 subject to 
\begin{multline*}
  f(g m m' e^{tH} n) 
  = 
  \sigma(m)^{-1} \delta(m') e^{-\lambda t} f(g)
\\
\text{for all
 $g \in G$, $m m' \in M \simeq SO(n) \times O(1)$, $t \in {\mathbb{R}}$, 
 $n \in N_+$.}
\end{multline*}
If $V$ is the representation of $SO(n)$
 on the exterior tensor space $\Exterior^i({\mathbb{C}}^{n})$
 ($2i \ne n$), 
 we use the notation $I_{\delta}(i, \lambda)$
 for $I_{\delta}(V, \lambda)$.  
Then the $SO(n)$-isomorphism
 on the exterior representations
$
  \Exterior^i({\mathbb{C}}^{n}) \simeq \Exterior^{n-i}({\mathbb{C}}^{n})
$
leads us to the following $G$-isomorphism:
\[
   I_{\delta}(i, \lambda) \simeq I_{\delta}(n-i, \lambda).  
\]
If $n$ is even and $n=2i$, 
the exterior representation $\Exterior^i({\mathbb{C}}^{n})$
 splits into two irreducible representations
 of $SO(n)$:
\[
{\Exterior}^{\frac n 2}({\mathbb{C}}^{n})
\simeq
{\Exterior}^{\frac n 2}({\mathbb{C}}^{n})_+
\oplus
{\Exterior}^{\frac n 2}({\mathbb{C}}^{n})_-
\]
with highest weights $(1,\cdots,1,1)$ and $(1,\cdots,1,-1)$, 
 respectively, 
 with respect to a fixed positive system for ${\mathfrak {so}}(n,{\mathbb{C}})$.  
Accordingly,
 we have a direct sum decomposition
 of the induced representation:
\begin{align}
   \operatorname{Ind}_{P}^{G}
  ({\Exterior}^{\frac n 2}({\mathbb{C}}^{n}) \otimes \delta, \lambda)
&=
  I_{\delta}({\Exterior}^{\frac n 2}({\mathbb{C}}^{n})_+, \lambda)
\oplus
   I_{\delta}({\Exterior}^{\frac n 2}({\mathbb{C}}^{n})_-, \lambda), 
\notag
\intertext{which we shall write as}
  I_{\delta}\left( \frac n 2, \lambda \right)
  &=
  I_{\delta}^{(+)} \left( \frac n 2, \lambda \right)
\oplus
   I_{\delta}^{(-)} \left( \frac n 2, \lambda \right) .  
\label{eqn:split}
\end{align}
Via the Harish-Chandra isomorphism, 
 the ${\mathfrak {Z}}({\mathfrak {g}})$-infinitesimal character
 of the trivial one-dimensional representation ${\bf{1}}$
 is given by 
\[
\rho= \left( \frac{n}{2}, \frac{n}{2}-1,\cdots,\frac{n}{2}-[\frac{n}{2}]
      \right)
\]
 in the standard coordinates
 of the Cartan subalgebra of
$
   {\mathfrak {g}}_{\mathbb{C}}
$
$
   =
$
$
   {\mathfrak {so}}(n+2,{\mathbb{C}})
$, 
 whereas that of $I_{\delta}(i, \lambda)$
 and $I_{\delta}^{(\pm)}(i, \lambda)$
 (when $n=2i$)
 is given by 
\begin{equation}
\label{eqn:Zginf}
  \left(
   \frac n 2, \frac n 2-1, \cdots, \frac n 2-i+1, 
   \widehat{\frac n 2 -i}, 
   \frac n 2-i-1, 
   \cdots, 
   \frac n 2 - [\frac n 2], 
   \lambda-\frac n 2
  \right).  
\end{equation}

For the group $G'=SO(n,1)$, 
 we shall use the notation $J_{\varepsilon}(j,\nu)$
 for the unnormalized parabolic induction
$
   \operatorname{Ind}_{P'}^{G'}
   (\Exterior^j({\mathbb{C}}^{n-1})\otimes \varepsilon,\nu)
$
 for $0 \le j \le n-1$, 
 $\varepsilon \in \{\pm\}$, 
 and $\nu \in {\mathbb{C}}$.

\subsection{Classification of symmetry breaking operators}
\label{subsec:allSBO}

Let $(G,G')=(SO(n+1,1), SO(n,1))$
 with $n \ge 3$.  
In this section we provide 
 a complete classification
 of symmetry breaking operators from 
 $I_{\delta}(i,\lambda)$ to $J_{\varepsilon}(j,\nu)$.  
The two recent articles \cite{KKP} and \cite{KS} gave 
 an explicit construction 
 and the classification
 of symmetry breaking operators in the following settings.  
\begin{enumerate}
\item[(1)]
$i=j=0$.  
The classification was accomplished in \cite{KS}.  
\item[(2)]
{\it{Differential}} symmetry breaking operators
 for $i$, $j$ general.  
The classification was accomplished in \cite[Thm.~2.8]{KKP}.  
\end{enumerate}

The proof of our general case 
 (Theorem \ref{thm:1.1} below)
 relies partially on the results 
 and techniques 
 that are developed in \cite{KKP, KS}.  
We note that the above literature treats the pair
 $(O(n+1,1), O(n,1))$, from 
which one can readily deduce the classification 
 for the pair
 $(SO(n+1,1), SO(n,1))$
 as we explained in \cite[Chap.~2, Sec.~5]{KKP}.

For the admissible smooth representations
 $\Pi=I_{\delta}(i,\lambda)$ of $G=SO(n+1,1)$
 and $\pi=J_{\varepsilon}(j,\nu)$ of $G'=SO(n,1)$, 
 we set
\[
   m(i,j)
   \equiv m(I_{\delta}(i,\lambda), J_{\varepsilon}(j,\nu))
   :=
   \dim \operatorname{Hom}_{G'}
    (I_{\delta}(i,\lambda)|_{G'}, J_{\varepsilon}(j,\nu)).  
\]
%

In order to give a closed formula
 of $m(i,j)$ as a function of $(\lambda, \nu,\delta, \varepsilon)$, 
 we introduce the following subsets of ${\mathbb{Z}}^2 \times \{\pm 1\}$:
\begin{align*}
  L:=&\{(-i,-j, (-1)^{i+j}): (i,j) \in {\mathbb{Z}}^2, 0 \le j \le i \}, 
\\
  L':=&\{(\lambda,\nu,\gamma) \in L: \nu \ne 0 \}.  
\end{align*}
To simplify the notation we also use will  use $\varepsilon$, $ \delta \in \{\pm \}$  and   $\varepsilon$, $ \delta \in \{\pm 1 \}$.

In the theorem below, we shall see
\begin{alignat*}{2}
m(i,j) \in &\{ 1,2,4 \} \qquad
&&\text{if $j=i-1$ or $i$}, 
\\
m(i,j) \in &\{ 0,1,2 \} \qquad
&&\text{if $j=i-2$ or $i+1$}, 
\\
m(i,j) =& 0 \qquad
&&\text{otherwise}.  
\end{alignat*}
Here is an explicit formula of the multiplicity
 for the restriction of nonunitary principal series representations
in this setting:

\begin{theorem}
\label{thm:1.1}
Suppose $n \ge 3$, $0 \le i\le [\frac n2]$, $0 \le j \le [\frac{n-1}2]$,
 $\delta$, $\varepsilon \in \{\pm \}\equiv \{ \pm 1\}$, 
 and $\lambda,\nu \in {\mathbb{C}}$.
Let $\Pi=I_{\delta}(i,\lambda)$
 and $\pi=J_{\varepsilon}(j,\nu)$
 be the admissible representations of $G=SO(n+1,1)$ and $G'=SO(n,1)$, respectively, as before.  
Then we have the following.  
\begin{enumerate}
\item[{\rm{(1)}}]
Suppose $j=i$.  
\begin{enumerate}
\item[{\rm{(a)}}]
Case $i=0$.
\begin{equation*}
m(0,0)
=
\begin{cases}
2
\qquad
&\text{if }
(\lambda, \nu, \delta\varepsilon) \in L, 
\\
1
&
\text{otherwise.}
\end{cases}
\end{equation*}
\item[{\rm{(b)}}]
Case $1 \le i < \frac  n2-1$.  
\begin{equation*}
m(i, i)
=
\begin{cases}
2
\qquad
&\text{if }
(\lambda, \nu, \delta \varepsilon) \in L' \cup \{(i,i,+)\},  
\\
1 
&
\text{otherwise}.  
\end{cases}
\end{equation*}
\item[{\rm{(c)}}]
Case $i= \frac n 2 -1$ ($n$: even).   
\begin{equation*}
m(\frac n 2 -1, \frac n 2 -1)
=
\begin{cases}
2
\qquad
&\text{if }
(\lambda, \nu, \delta\varepsilon) \in L' \cup \{(i,i,+)\}\cup \{(i,i+1,-)\},   
\\
1
&
\text{otherwise}.  
\end{cases}
\end{equation*}
\item[{\rm{(d)}}]
Case $i= \frac {n-1} 2$ ($n$: odd).   
\begin{equation*}
m(\frac {n-1} 2,\frac {n-1} 2)
=
\begin{cases}
4
\qquad
&\text{if }
(\lambda, \nu, \delta\varepsilon) \in L' \cup \{(i,i,+)\}, 
\\
2
&
\text{otherwise}.  
\end{cases}
\end{equation*}
\end{enumerate}

\item[{\rm{(2)}}]
Suppose $j=i-1$.  
\begin{enumerate}
\item[{\rm{(a)}}]
Case $1 \le i < \frac{n-1}{2}$.
\begin{equation*}
m(i,i-1)
=
\begin{cases}
2
\qquad
&\text{if }
(\lambda, \nu, \delta\varepsilon) \in L' \cup \{(n-i,n-i,+)\}, 
\\
1
&
\text{otherwise.}
\end{cases}
\end{equation*}
\item[{\rm{(b)}}]
Case $i = \frac{n-1}{2}$ ($n$: odd).  
\begin{equation*}
m(\frac{n-1}{2},\frac{n-3}{2})
=
\begin{cases}
2
\qquad
&\text{if }
(\lambda, \nu, \delta \varepsilon) \in L' \cup \{(n-i,n-i,+)\}
                                       \cup \{(i,i+1,-)\} ,  
\\
1 
&
\text{otherwise}.  
\end{cases}
\end{equation*}
\item[{\rm{(c)}}]
Case $i= \frac n 2$ ($n$: even).   
\begin{equation*}
m(\frac n 2,\frac n 2-1)
=
\begin{cases}
4
\qquad
&\text{if }
(\lambda, \nu, \delta\varepsilon) \in L' \cup \{(n-i,n-i,+)\},   
\\
2
&
\text{otherwise}.  
\end{cases}
\end{equation*}
\end{enumerate}

\item[{\rm{(3)}}]
Suppose $j=i-2$.  
\begin{enumerate}
\item[{\rm{(a)}}]
Case $2 \le i < \frac{n}{2}$.
\begin{equation*}
m(i,i-2)
=
\begin{cases}
1
\qquad
&\text{if }
(\lambda, \nu, \delta\varepsilon) = (n-i,n-i+1,-), 
\\
0
&
\text{otherwise.}
\end{cases}
\end{equation*}
\item[{\rm{(b)}}]
Case $i = \frac{n}{2}$ ($n$: even).  
\begin{equation*}
m(\frac{n}{2},\frac{n}{2}-2)
=
\begin{cases}
2
\qquad
&\text{if }
(\lambda, \nu, \delta \varepsilon) =(\frac{n}{2},\frac{n}{2}+1,-),  
\\
0 
&
\text{otherwise}.  
\end{cases}
\end{equation*}
\end{enumerate}

\item[{\rm{(4)}}]
Suppose $j=i+1$.  
\begin{enumerate}
\item[{\rm{(a)}}]
Case $i=0$.
\begin{equation*}
m(0,1)
=
\begin{cases}
1
\qquad
&\text{if }
\lambda \in - {\mathbb{N}}, \nu=1, \text{ and }\delta\varepsilon = (-1)^{\lambda+1}, 
\\
0
&
\text{otherwise.}
\end{cases}
\end{equation*}
\item[{\rm{(b)}}]
Case $1 \le i < \frac{n-3}{2}$.  
\begin{equation*}
m(i,i+1)
=
\begin{cases}
1
\qquad
&\text{if }
(\lambda, \nu, \delta \varepsilon) =(i,i+1,-),  
\\
0 
&
\text{otherwise}.  
\end{cases}
\end{equation*}
\item[{\rm{(c)}}]
Case $i=\frac{n-3}{2}$ ($n$: odd).  
\begin{equation*}
m(\frac{n-3}{2}, \frac{n-1}{2})
=
\begin{cases}
2
\qquad
&\text{if }
(\lambda, \nu, \delta \varepsilon) =(\frac{n-3}{2},\frac{n-1}{2},-),  
\\
0 
&
\text{otherwise}.  
\end{cases}
\end{equation*}
\end{enumerate}

\item[{\rm{(5)}}]
Suppose $j \not\in \{i-2, i-1, i, i+1\}$.
Then $m(i,j)=0$ for all $\lambda, \nu, \delta, \varepsilon$.
\end{enumerate}

\end{theorem}

The construction of nontrivial symmetry breaking operators
 is proved by generalizing the techniques developed in \cite{KS} 
 in the scalar case to the matrix-valued case
 for representations induced from finite-dimensional representations of $M$. 
The proof for the exhaustion of (continuous) symmetry breaking operators
 is built on the classification of differential symmetry breaking operators
 which was given in \cite[Thm.~2.8]{KKP}.

\begin{remark}
[multiplicity-one property]
\label{rem:mult2}
In \cite{SunZhu}
 it is proved that 
$$
   \dim_{\mathbb{C}}\operatorname{Hom}_{G'}(\Pi|_{G'},\pi) \le 1
$$
 for any irreducible admissible smooth representations
 $\Pi$ and $\pi$ of $G=SO(n+1,1)$ and $G'=SO(n,1)$, 
respectively.  
Thus Theorem \ref{thm:1.1} fits well 
 with their multiplicity-free results 
 for 
  $\lambda, \nu \in {\mathbb{C}}\setminus {\mathbb{Z}}$, 
 where $I_{\delta}(i,\lambda)$ and $J_{\varepsilon}(j,\nu)$
 are irreducible admissible representations
 of $G$ and $G'$, 
 respectively,
 except for the cases $n=2 i$ or $n= 2j+1$.  
We note that, 
 in addition to the subgroup $G'=SO(n,1)$, 
 the Lorentz group contains two other subgroups
 of index two, 
 that is, 
 $O^+(n,1)$ (containing orthochronous reflections)
 and  $O^-(n,1)$ (containing anti-orthochronous reflections)
 with terminology
 in relativistic space-time for $n=3$.  
Our method gives also the multiplicity formula for such pairs,
and it turns out that an analogous multiplicity-one statement
 fails 
 if we replace $(G,G')=(SO(n+1,1),SO(n,1))$
 by $(O^-(n+1,1),O^-(n,1))$.  
In fact, 
 the multiplicity $m(\Pi, \pi)$ may equal 2
 for irreducible representations $\Pi$ and $\pi$
 of $O^-(n+1,1)$ and $O^-(n,1)$, 
 respectively.  
\end{remark}

\bigskip
\section{ \bf Main Results: Symmetry breaking for representations of rank one orthogonal groups}
\label{sec:mainresults}

The main result in this section is a theorem about multiplicities for irreducible representations with trivial infinitesimal character $\rho$, 
 namely, 
those representations
 that have the same ${\mathfrak{Z}}({\mathfrak {g}})$-infinitesimal 
 character
 with the trivial one-dimensional representation.  
We first state the result  using the Langlands parameters of the irreducible representations \cite{BW, L}. In the second part we introduce $\theta $-stable parabolic pairs $\bq, L$ and parametrize the representations
 by one-dimensional representations of $L$
 following \cite{KV, KMemoirs92, VZ}. We then state again the theorem in this formalism.

\subsection{Irreducible representations with infinitesimal \\  character $\rho$}\label{subsec:isSO}

In this section
 we give a description
 of all irreducible admissible representations
 of $G=SO(n+1,1)$
 with trivial infinitesimal character $\rho$.  
Another description will be given in Section \ref{subsec:main2}.

By \eqref{eqn:Zginf}, 
 $I_{\delta}(i, \lambda)$ has the ${\mathfrak {Z}}({\mathfrak {g}})$-infinitesimal character 
 $\rho$ 
 if and only if $\lambda=i$ or $\lambda=n-i$.  
We identify the maximal compact subgroup $K$ of $G$
 with $O(n+1)$
 via the isomorphism \eqref{eqn:K}.  
In what follows,
 we use the notation of \cite[Chap.~2, Sect.~3]{KKP}
 by adapting it to $SO(n+1,1)$
 instead of $O(n+1,1)$.  
For $0 \le i \le n$, 
 we denote by $I_{\delta}(i)^{\flat}$ and $I_{\delta}(i)^{\sharp}$
 the unique irreducible subquotients
 of $I_{\delta}(i, i)$
 containing the irreducible representations
 $\Exterior^i({\mathbb{C}}^{n+1})\otimes \delta$
 and $\Exterior^{i+1}({\mathbb{C}}^{n+1})\otimes (-\delta)$
 of $O(n+1) \simeq K$, 
 respectively.  
In the case $n=2i$, 
 the $SO(n+1,1)$-modules
 $I_{\delta}^{(\pm)}(\frac n 2, \frac n2)$
 are irreducible for $\delta= \pm$, 
 and we have the following isomorphisms:
\begin{equation}
\label{eqn:1616116}
   I_{\delta}\left( \frac n 2 \right)^{\flat}
   \simeq 
   I_{\delta}\left( \frac n 2 \right)^{\sharp}
   \simeq
   I_{\delta}^{(+)}\left( \frac n 2, \frac n2 \right)
   \simeq 
   I_{\delta}^{(-)}\left( \frac n 2, \frac n2 \right)
\qquad
 \text{for $\delta=\pm$, }
\end{equation}
as representations
 of $SO(n+1,1)$.

Then we have $G$-isomorphisms:
\begin{equation}
\label{eqn:LNS218}
  I_{\delta}(i)^{\sharp} \simeq I_{-\delta}(i+1)^{\flat}
\quad
  \text{for $0 \le i \le n$ and $\delta \in \{\pm \}$}.  
\end{equation}
For $0 \le \ell \le n+1$ and $\delta \in \{\pm \}$, 
 we set 
\begin{equation}
\label{eqn:Pild}
  \Pi_{\ell,\delta}
  :=
  \begin{cases}
  I_{\delta}(\ell)^{\flat} \quad&(0 \le \ell \le n), 
  \\
  I_{-\delta}(\ell-1)^{\sharp} \quad&(1 \le \ell \le n+1).  
  \end{cases}
\end{equation}
In view of \eqref{eqn:LNS218}, 
 $\Pi_{\ell,\delta}$ is well-defined.

For $0 \le i \le n$ with $n \ne 2i$ and $\delta \in \{ \pm \}$, 
 we have a nonsplitting exact sequence of $G$-modules:
\[
0 \to \Pi_{i,\delta} \to I_{\delta}(i,i) \to \Pi_{i+1,-\delta}
\to 0.  
\]
As we mentioned in \eqref{eqn:split}, 
 $I_{\delta}(\frac n 2, \frac n2)=\Pi_{\frac n2, \delta}\oplus \Pi_{\frac n2+1, -\delta}$
 when $n=2i$.  
 
 \medskip
The properties
 of irreducible representations
 $\Pi_{\ell,\delta}$
 ($0 \le \ell \le n+1$, $\delta = \pm$)
can be summarized as follows \cite{BW, KKP}.

\begin{prop}
\label{prop:161648}
Let $G:=SO(n+1,1)$
 with $n \ge 1$.  
\begin{enumerate}
\item[{\rm{(1)}}]
$\Pi_{\ell,\delta} \simeq \Pi_{n+1-\ell,-\delta}$
 as $G$-modules
 for all $0 \le \ell \le n+1$ and $\delta = \pm $.  
\item[{\rm{(2)}}]
Irreducible admissible representations of moderate growth 
 with ${\mathfrak {Z}}({\mathfrak {g}})$-infinitesimal character $\rho$
 are classified as 
\begin{alignat*}{2}
&
 \{ 
 \Pi_{\ell,\delta}: 0 \le \ell \le \frac{n-1}{2}, \delta = \pm 
 \}
 \cup
 \{ 
 \Pi_{\frac{n+1}{2},+}
 \}
 \qquad
&&
 \text{if $n$ is odd}, 
\\
&
 \{ 
 \Pi_{\ell,\delta}: 0 \le \ell \le \frac{n}{2}, \delta = \pm 
 \}
 \qquad
&&
 \text{if $n$ is even}.  
\end{alignat*}

\item[{\rm{(3)}}] Every $\Pi_{\ell,\delta}$ is unitarizable.  
\end{enumerate}
By abuse of notation, 
 we use the same symbol $\Pi_{l,\delta}$
 to denote the unitarization.  

\begin{enumerate}
\item[{\rm{(4)}}] For $n$ odd, $\Pi_{\frac{n+1}{2},+}$ is
 a discrete series.  
For $n$ even, 
 $\Pi_{\frac{n}{2},\pm}$ are 
 tempered representations.  
All the other representations in the list (2) are nontempered.  
\item[{\rm{(5)}}]
For $n$ even, 
 the center of $G$ acts nontrivially on $\Pi_{\ell,\delta}$
 if and only if $\delta=(-1)^{\ell+1}$.  
For $n$ odd, 
 the center of $G$ is trivial,  
 and thus acts trivially on $\Pi_{\ell,\delta}$
 for any $\ell$ and $\delta$.  
\end{enumerate}
\end{prop}

For the subgroup $G'=SO(n,1)$, 
 we shall use similar notation $\pi_{j,\varepsilon}$
 for the subrepresentations
 of $J_{\varepsilon}(j,j)$
 (or the quotients of $J_{\varepsilon}(j-1,j-1)$).

In view of Proposition \ref{prop:161648}, 
 in particular, the $G$-isomorphism
 $\Pi_{\frac{n+1}{2}, +} \simeq \Pi_{\frac{n+1}{2}, -}$ for $n$ odd 
 and the $G'$-isomorphism
 $\pi_{\frac{n}{2}, +} \simeq \pi_{\frac{n}{2}, -}$ for $n$ even, 
 we shall use the following convention: 
\begin{enumerate}
\label{eqn:halfsgn}
\item  {if $n+1=2i$; we identify $\delta =+ $ and $\delta = -$ 
\item  if $n=2j$ we identify $\varepsilon =+ $ and $\varepsilon =-$ }
\end{enumerate} 
in statements and theorems  about  representations
$\Pi_{i,\delta}$
  and   $\pi_{j,\varepsilon}$ with indices ($0 \le i \le [\frac{n+1}{2}]$) and 
($0 \le j \le [\frac{n}{2}]$).  

\subsection{Formulation I of the main theorem} 
\label{subsec:formulationI}

As we saw in Proposition \ref{prop:161648}, 
 all the  representations $\Pi_{i,\delta}$ are unitarizable,  
 but  the restriction of  $\Pi_{i,\delta}$ to the subgroup $G'$  does not decompose into a direct sum of irreducible representations
 \cite{KInvent98}.  
Hence to obtain information about the restriction we consider 
 $G'$-intertwining operators
 ({\it{symmetry breaking operators}})
 for smooth admissible representations:
\begin{equation}
\label{eqn:Hom}
 \operatorname{Hom}_{G'} ({\Pi_{i,\delta}}|_{G'}, \pi_{j,\varepsilon}) .
\end{equation}

By using the classification of all symmetry breaking operators for principal series representations
(Theorem~\ref{thm:1.1}) and by analyzing their restrictions to the subquotients
of principal series representations, 
 we can determine the dimension of the space \eqref{eqn:Hom},
 and, in particular, we obtain a necessary and sufficient condition for this space
 to be nonzero.  
Here is a statement.

\begin{theorem}
\label{thm:1616112}
Let $(G,G')=(SO(n+1,1),SO(n,1))$.  
Suppose $0 \le i \le [\frac {n+1}{2}]$, 
 $0 \le j \le [\frac {n}{2}]$, 
 and $\delta$, $\varepsilon= \pm$
 with the convention \eqref{eqn:halfsgn}.  
Then 
\[
  \dim_{\mathbb{C}} \operatorname{Hom}_{G'}
   (\Pi_{i,\delta}|_{G'}, \pi_{j,\varepsilon})
  =
  \begin{cases}
   1
   \qquad
   & \text{if $\delta=\varepsilon$ and $j \in \{ i-1,i \}$,}
\\
   0
   \quad
   & \text{otherwise}.  
  \end{cases}
\]
\end{theorem}
Theorem \ref{thm:1616112} can be rephrased as follows. 
\begin{theorem} 
\label{symbreaking/rho}
Suppose $0 \le i \le [\frac{n+1}{2}]$,
  $0 \le j \le [\frac{n}{2}]$, 
 and $\delta, \varepsilon= \pm$.  
Then the following three conditions
on the quadruple
 $(i,j,\delta,\varepsilon)$
 are equivalent.  
\begin{enumerate}
\item[{\rm{(i)}}]
\[
   \operatorname{Hom}_{G'} (\Pi_{i,\delta}|_{G'},\pi_{j,\varepsilon})\ne\{0\}. 
\]
\item[{\rm{(ii)}}]
\[
   \dim_{\mathbb{C}}\operatorname{Hom}_{G'} (\Pi_{i,\delta}|_{G'},\pi_{j,\varepsilon})  = 1.  
\]
\item[{\rm{(iii)}}]
There is an arrow connecting the representations in the following tables
 with $\delta=\varepsilon$. 
(For simplicity, 
 we omit the subscripts $\delta$ and $\varepsilon$
 in the tables below.)
\end{enumerate}
 \medskip
In (iii), the convention \eqref{eqn:halfsgn} is applied
 to the cases $j=m$ when $n=2m$ (see Table \ref{default})
 and $i=m+1$ when $n=2m+1$ (see Table \ref{default2}), 
 where $\pi_{m,+} \simeq \pi_{m,-}$
 and $\Pi_{m+1,+} \simeq \Pi_{m+1,-}$ hold, 
 respectively.
\begin{table}[htp]
\caption{Symmetry breaking for $(SO(2m+1,1), SO(2m,1))$ }

\begin{center}
\begin{tabular}{cccccccccc}
$\Pi_0$& &$\Pi_1$& &\dots &\dots \ \ \ \ \ & \dots & $\Pi_{m-1} $& & $\Pi_{m}$ \\
$\downarrow$ &$\swarrow$& $\downarrow $& $\swarrow$ & \dots  & & & $ \downarrow $&  $\swarrow $  &  $\downarrow$ \\
$\pi_0$& &$\pi_1$& &\dots  &\dots \ \ \ \ \ &\dots & $\pi_{m-1}$ & & $\pi_{m}$ 
\end{tabular}
\end{center}
\label{default}
\end{table}

 \medskip
\begin{table}[htp]
\caption{Symmetry breaking for $(SO(2m+2,1),SO(2m+1,1))$ }
\begin{center}
\begin{tabular}{ccccccccccc}
$\Pi_0$& &$\Pi_1$& &\dots &\dots \ \ \ \ \ & \dots & $\Pi_{m-1} $& & $\Pi_{m}$ & $\Pi_{m+1}$\\
$\downarrow$ &$\swarrow$& $\downarrow $& $\swarrow$ & \dots  & & & $ \downarrow $& $\swarrow $ & $\downarrow$ & $\swarrow$\\
$\pi_0$& &$\pi_1$& &\dots  &\dots \ \ \ \ \ &\dots & $\pi_{m-1}$ & & $\pi_{m}$ 
\end{tabular}
\end{center}
\label{default2}
\end{table}
\end{theorem}

We note that the equivalence (i)  $\Leftrightarrow$ (ii) 
in Theorem \ref{symbreaking/rho} 
 could be derived also from 
 the general theory \cite{SunZhu} because $\Pi_{i,\delta}$ and $\pi_{j,\varepsilon}$
 are irreducible.

\subsection{Formulation II of the main theorem} 
\label{subsec:main2}
We have described the irreducible representations $\Pi_{\ell,\delta}$
 as subquotients
 of principal series representations
 in Section \ref{subsec:isSO}.  
The next proposition provides another characterization
 of the same representations $\Pi_{\ell,\delta}$.  

\begin{prop}
The irreducible representations $\Pi_{\ell,\delta}$ in Proposition \ref{prop:161648} (3)
 are the Casselman--Wallach 
globalization
 of the irreducible, 
 unitarizable $({\mathfrak {g}}, K)$-modules
 with nonzero $({\mathfrak {g}}, {\mathfrak {k}})$-cohomologies.  
\end{prop}

These $({\mathfrak {g}}, K)$-modules
 can be described by using the Zuckerman derived functor modules.  
For this,
 let us introduce some notation.  
For $0 \le i \le [\frac {n+1}{2}]$, 
 we consider a $\theta$-stable parabolic subalgebra
 ${\mathfrak {q}}_i$ 
 of ${\mathfrak {g}}_{\mathbb{C}} = {\mathfrak {so}}(n+1,1) \otimes_{\mathbb{R}}{\mathbb{C}} \simeq {\mathfrak {so}}(n+2, {\mathbb{C}})$ 
with the (real) Levi subgroup
\[
   L_i := N_G({\mathfrak {q}}_i)
       \simeq SO(2)^i \times SO(n+1-2i, 1).  
\]
We note
 that $L_i$ meets all the connected components
 of $G=SO(n+1,1)$.  
For the trivial one-dimensional representation 
 ${\bf{1}}$ 
 of the first factor $SO(2)^i$
 and a one-dimensional representation $\chi$
 of the last factor $SO(n+1-2i, 1)$, 
 we define a $({\mathfrak {g}}, K)$-module
\[
   A_{{\mathfrak {q}}_i}(\chi)
   :=
   {\mathcal{R}}_{{\mathfrak {q}}_i}^{S_i}({\bf{1}} \boxtimes \chi)
\]
  as the cohomological parabolic induction from 
 the one-dimensional representation 
 ${\bf{1}} \boxtimes \chi$
 of $L$.  
We adopt a \lq{$\rho$-shift}\rq\ of the cohomological parabolic induction 
 in a way 
 that $A_{{\mathfrak {q}}_i}(\chi)$ has 
 the infinitesimal character $\rho$
 if $d \chi =0$.  
(The $({\mathfrak {g}}, K)$-module 
 $A_{{\mathfrak {q}}_i}(0)$
 in the notation of Vogan--Zuckerman \cite{KO}
 corresponds to $A_{{\mathfrak {q}}_i}({\bf{1}})$
 in our notation.)  
There are two characters 
 $\chi$ of $SO(k,1)$
 ($k \ge 1$)
 such that $d\chi=0$.  
We write $\chi_{k,1}^+$ for the trivial one,
 and $\chi_{k,1}^-$ for the nontrivial one.  
Then we have
\begin{prop}
[{\cite{KV, KMemoirs92}}]
\label{prop:161655}
Suppose $0 \le i \le [\frac {n+1}{2}]$.  
For $\varepsilon =\pm$, 
 
\begin{equation*}
  (\Pi_{i,\varepsilon})_K \simeq  A_{{\mathfrak {q}}_i}(\chi_{n+1-2i,1}^{\varepsilon}).  
\end{equation*}
\end{prop}

\begin{remark}
{\rm{
We may regard $\chi_{0,1}^- \simeq \chi_{0,1}^+$
 for the representation of $SO(0,1)=\{1\}$.  
When $n$ is odd
 and $i = \frac{n+1}{2}$, 
$L \simeq SO(2)^{\frac{n+1}{2}} \times SO(0,1)$.  
This matches the $({\mathfrak {g}}, K)$-isomorphism:
 $\Pi_{\frac{n+1}{2},+} \simeq \Pi_{\frac{n+1}{2},-}$
 (see Proposition \ref{prop:161648} (1)).  
}}
\end{remark}

\medskip
\begin{example}
[see Proposition \ref{prop:161648}]
~~~
{\rm{
\begin{itemize}
\item[{\rm{(1)}}] $\Pi_{0,\varepsilon}$ is one-dimensional.  

\item[{\rm{(2)}}] If $n = 2m$ then $\Pi_{m,+}$ and $\Pi_{m,-}$ are
 the inequivalent tempered principal series representations
 of $SO(2m+1,1)$ with infinitesimal character $\rho$.  

\item[{\rm{(3)}}]  If $n= 2m-1$ then $\Pi_{m,+} \simeq \Pi_{m,-}$
 is the unique discrete series representation
 of $SO(2m,1)$ with infinitesimal character $\rho$.
\end{itemize}
}}
\end{example}

\medskip 
\subsection{$\theta$-stable parameter of $\Pi_{i,\delta}$}
\label{subsec:III.4}
Suppose $0 \le i \le [\frac{n+1}{2}]$. 
Let $\Sigma_i^+$ be the set of positive roots
 corresponding to the nilpotent radical
 of the $\theta$-stable parabolic subalgebra ${\mathfrak{q}}_i$
 and define 
\[\rho_i= \frac 1 2 \sum _{\alpha \in \Sigma_i^+} \alpha .\]
Via the standard basis
 of the fundamental Cartan subalgebra, 
 we have
\begin{alignat*}{2}
\rho_i&=( m, m-1,  \dots , m-i+1, 0 \dots 0 )
\quad
&&\mbox{if $G=SO(2m+1,1)$, } 
\\
 \rho_i&=  (m- \frac 1 2, m- \frac 3 2,  \dots , m-i - \frac 1 2, 0,\dots  0 ) 
\quad
&&\mbox{if $G=SO(2m,1)$.} 
\end{alignat*}
To make our notation consistent with the Harish-Chandra parameter for discrete series representations for $SO(2m,1) $ 
 we define the {\it{$\theta$-stable parameters}}  of the  cohomologically induced representation 
 $(\Pi_{i,\delta})_K\simeq A_{{\mathfrak {q}}_i}(\chi_{n-2i+1,1}^{\delta})$
 as follows.  

\medskip
\begin{definition}
Suppose $0 \le i \le m$ and $\delta \in \{\pm\}$.  

\begin{itemize}
\item[{\rm{(1)}}]
The {\it{$\theta$-stable parameter}}
 of the irreducible representation $\Pi_{i,\delta}$  of  \linebreak $SO(2m+1,1)$
 is 
\[
   (m, m-1,  \dots , m-i+1  \  || \ \chi_{2m-2i+1,1}^{\delta}), 
\]
where $\chi_{2m-2i+1,1}^{\delta}$ is
 the one-dimensional representation of  \linebreak $SO(2m-2i+1,1)$. 
\item[{\rm{(2)}}]
The {\it{$\theta$-stable parameter}} 
 of the irreducible representation $\Pi_{i,\delta}$ of $SO(2m,1)$ is 
\[
(m- \frac 1 2, m- \frac 3 2,  \dots , m-i +\frac 1 2 \  || \ \chi_{2m-2i,1}^
{\delta}), 
\]
where $\chi_{2m-2i,1}^{\delta}$ is the one-dimensional representation of \linebreak $SO(2m-2i,1)$.  
\end{itemize}
\end{definition}

\medskip \noindent
We use the same convention
 for the representations $\pi_{j,\varepsilon}$
 of $G'$.

\medskip

Theorem \ref{symbreaking/rho}  can now be restated in a formulation resembling the classical branching law for finite-dimensional representations. 
We connect the parameter
 by an arrow $\Downarrow$ pointing towards the parameter
 of the representation of the smaller group.


\begin{theorem} 
Suppose that $(G,G')=(SO(n+1,1), SO(n,1))$.  
Let $\Pi$ and $\pi$ be irreducible admissible representations
 of moderate growth of $G$ and $G'$
 with ${\mathfrak{Z}}({\mathfrak{g}})$-infinitesimal character $\rho$, 
 respectively.  
\begin{enumerate}
\item[{\rm{(1)}}]   
Suppose $n=2m$.  
Then
\[
  \dim_{\mathbb{C}} 
  \operatorname{Hom}_{G'} (\Pi|_{G'},\pi)  
  = 1 
\]
if and only if the  $\theta$-stable parameters of $\Pi$ and $\pi$
 satisfy one of the following conditions:
\\
 $(\Pi, \pi)=(\Pi_{i,\delta}, \pi_{i,\varepsilon})$
 for $0 \leq i \leq m$ with $\delta=\varepsilon \in \{\pm \}$
 (the convention \eqref{eqn:halfsgn} is applied to $\varepsilon$
 when $i=m$): 
\begin{eqnarray*}
& ( m ,  m-1, \dots   m+1-i \ || \  \chi_{2m+1-2i ,1}^{\delta} )  &    \\
&         \Downarrow      &  \\
 &  (m-\frac 12, m- \frac 32,  \dots , m +\frac 12 -i   \  || \ \chi_{2m-2i,1}^{\varepsilon} )&
\end{eqnarray*}
or $(\Pi, \pi)=(\Pi_{i,\delta}, \pi_{i-1,\varepsilon})$ for $0<i \leq m$
 with $\delta =  \varepsilon \in \{\pm\}$: 
\begin{eqnarray*}
&( m ,  m-1, \dots   m+1-i \ || \  \chi_{2m+1-2i,1}^{\delta})  &   \\
&  \Downarrow &  \\
  & (m-\frac 12, m-\frac32 , \dots ,m+\frac 32-i \ || \ \chi_{2m+2-2i,1}^{\varepsilon} ).   &
\end{eqnarray*}
\item[{\rm{(2)}}]
Suppose $n=2m+1$.  
Then 
\[
     \dim_{\mathbb{C}} \operatorname{Hom}_{G'} (\Pi|_{G'},\pi) =1 
\]
if and only if the $\theta$-stable parameters of $\Pi$ and $\pi$ satisfy 
 one of the following conditions:
\\
 $(\Pi, \pi)=(\Pi_{i,\delta}, \pi_{i,\varepsilon})$
 for  $ 0 \leq i < m+1$
 with $\delta$, $\varepsilon \in \{\pm\}$: 
 \begin{eqnarray*}
&  (m+\frac 1 2, m- \frac 1 2,  \dots , m+ \frac 3 2-i \  || \ \chi_{2m+2-2i,1}^{\delta})  &   \\
& \Downarrow &  \\
&   (m,m-1,\dots , m+1-i \ || \ \chi_{2m+1-2i,1}^{\varepsilon}) 
\end{eqnarray*}
or $(\Pi, \pi)=(\Pi_{i,\delta}, \pi_{i-1,\varepsilon})$
 for $0 <i \leq m+1$ with $\delta= \varepsilon \in \{\pm\}$
 (the convention \eqref{eqn:halfsgn} is applied to $\delta$ when $i=m+1$):
\begin{eqnarray*}
 &  (m+ \frac 1 2, m- \frac 1 2,  \dots , m + \frac 3 2-i \  || \ \chi_{2m+2-2i,1}^{\delta}) &  \\
& \Downarrow  &  \\
& (m,m-1, \dots , m+2-i \ || \ \chi_{2m+3 -2i,1}^{\varepsilon}). &
\end{eqnarray*}
\end{enumerate}

\end{theorem}

\medskip

\noindent
\begin{remark}
The first case represents the vertical arrows and the second case represents the slanted arrow in Theorem \ref{symbreaking/rho} (iii).
\end{remark}
\medskip

\bigskip

\section{ \bf Symmetry breaking and the Gross--Prasad conjectures}
\label{sec:GrossPrasad}

In 2000 B.~Gross and N.~Wallach \cite{GW} showed
 that the restriction of {\it{small}} discrete series representations
 of $G=SO(2p+1,2q)$ to $G'=SO(2p,2q)$ 
satisfies the Gross--Prasad conjectures
 \cite{GP}. 
In that case,
 both the groups $G$ and $G'$ admit discrete series representations.  
On the other hand,
 for the pair $(G,G')=(SO(n+1,1),SO(n,1))$, 
 only one of $G$ or $G'$ admits discrete series representations.  
We sketch  here  a proof that  our theorem \ref{thm:1616112}
confirms the Gross--Prasad conjectures
 also for {\bf tempered} representations with infinitesimal character 
$\rho$. 

In our formulation and the exposition
 we rely on the original article
 by B.~Gross and D.~Prasad \cite{GP} and also on \cite{GGP}.

\medskip
The following diagram recalls our results
 in the previous sections 
 about symmetry breaking operators for tempered representations with infinitesimal character $\rho $ of the groups 
 $SO(n+1,1)$ for $n=2m$, $2m$, and $2m-1$. 
We denote the corresponding representations by $\Pi$, $\pi$ and $\varpi$, 
 respectively, 
 using the subscripts defined in Section \ref{sec:mainresults}.
For $n=2m$, 
 we simply write $\pi_m$ for $\pi_{m,\varepsilon}$
 because $\pi_{m,+} \simeq \pi_{m,-}$
 as $SO(2m,1)$-modules.  
\medskip
\begin{table}[h] 
\caption{Symmetry breaking for $SO(2m+1,1) \supset SO(2m,1) \supset SO(2m-1,1)$}\begin{center}
\begin{tabular}{ccccc}
$\Pi_{m,-}$&                  &             &                      &   $\varpi_{m-1,-}$ \\
             &$\searrow$     &             &  $\nearrow$ &                                 \\
            &                         &$ \pi_m $ &                   &                                 \\
             &$\nearrow $      &           & $\searrow$  &                                   \\
 $\Pi_{m,+}$&                  &            &                     & $\varpi_{m-1,+}$  
\end{tabular}
\end{center}
\label{default3}
\end{table}

We will in the following only consider representations which are nontrivial
 on the center {{(see Proposition \ref{prop:161648} (5))}}
 and thus are genuine representations of the orthogonal groups.
So we are considering in our discussion of the Gross--Prasad conjectures only this part of the diagram.

\bigskip
\begin{table}[h] 
\begin{center}
\begin{tabular}{ccccc}
{$\Pi_{m,(-1)^{m+1}}$}& $ \rightarrow   $   & $ \pi_m     $      &    $   \rightarrow      $         &   {$\varpi_{m-1,(-1)^m}$} \\
 &   &   &  &
\end{tabular}
\end{center}
\end{table}

\noindent
The other remaining cases can be handled by using the same ideas.

\medskip
 We first sketch  the results about Vogan packets for special orthogonal groups.
The Vogan $L$-packet 
 is the disjoint union of Langlands $L$-packet over {  \it pure } inner forms.  We refer to \cite{ABV} and \cite{V} for general information
 about Vogan packets and to \cite{GGP}
 for details  for special orthogonal groups.

\medskip
Consider the complexification $SO(n+1,\bC)$ of a special orthogonal group 
$SO(n,1)$ and let $T_\bC \subset SO(n+1,\bC)$ be the complexification of a fundamental Cartan subgroup $T$ of $SO(n+1,1)$. 

\medskip
\subsection {Vogan packets of discrete series representations with infinitesimal character $\rho $ of odd special orthogonal groups} 

We begin with the case
 $n=2m-1$.  
In this case $SO(n+1,1)=SO(2m,1)$
 has discrete series representations.  
We fix a set of positive roots 
$
   \Delta^+ \subset {\mathfrak{t}}_{\mathbb{C}}^{\ast}
$
 for the root system 
 $\Delta({\mathfrak{so}}(2m+1,\bC),{\mathfrak{t}}_{\mathbb{C}})$,
 and denote by $\rho$ half the sum of positive roots as before. 
For $l+k=2m+1$, 
 we call a real form $SO(l,k)$ {\it{pure}}
 if $l$ is even.
The Vogan packet containing the discrete series representation $\pi_m$ is the disjoint union of discrete series representations with infinitesimal character $\rho$ of the pure inner forms. 
The cardinality of this packet is 
\[
   2^m=\sum_{\substack{0 \le l \le 2m \\ l:\text{even}}} 
({m \atop{\frac l 2}}).  
\]
There exists a finite group 
$
  {\mathcal  A}_2 \simeq (\bZ_2)^m$  whose characters parametrize the representations in the Vogan packet. For the discrete series representation with parameter $\chi \in \widehat{{\mathcal{A}}_2}$
 we write $\pi(\chi). $ For more details
 see \cite{GR} or \cite{V}. 
We write $VP(\pi_m)$ for the Vogan packet containing $\pi_m $.

\medskip
\begin{example}
{\rm{
\begin{itemize}
\item[{\rm{(1)}}]
The trivial representation  ${\bf{1}}$ of  $SO(0,2m+1)$ is in $VP(\pi_m)$.

\item[{\rm{(2)}}] We can define similarly a Vogan packet containing $(SO(1,2m),\pi_m)$.  
\end{itemize}
}}
\end{example}

\medskip
By abuse of notation we may also consider $\pi_m$ as a discrete series representation of $SO(1,2m)$, but  the pairs $(SO(1,2m) ,\pi_m) $ and $(SO(2m,1),\pi_m)$ are not in the same Vogan packet.

\begin{remark}
{\rm{
Analogous results hold
 for the infinitesimal character $\lambda+\rho$ where $\lambda$ is 
the highest weight of a finite-dimensional representation.
}}
\end{remark}

\bigskip

\subsection{ A Vogan packet of tempered induced representations  with infinitesimal character $\rho$}
\label{subsec:IV.2}

Next we consider the case $n=2m$.  
Then $SO(n+1,1)=SO(2m+1,1)$ has no discrete series representation
 and we consider instead a Vogan packet of tempered representations
 with infinitesimal character $\rho$ which contains the pair $ (SO(2m+1,1) , \Pi_{m,\delta}) $ with $\delta= (-1)^{m+1}$.

To simplify the notation we assume  in this subsection that  $\delta =(-1)^{m+1}$. 
Recall that $\Pi_{m,\delta}$ denotes the irreducible representation 
$
   I_{\delta}(\Exterior^{m}({\mathbb{C}}^{2m})_+, m)
$
$
   \simeq
$
$
   I_{\delta}(\Exterior^{m}({\mathbb{C}}^{2m})_-, m)
$ 
 which is the smooth representation
 of a unitarily induced principal series representation from the maximal parabolic subgroup. 
Its Langlands/Vogan parameter factors through the Levi subgroup of a maximal parabolic subgroup of the Langlands dual group  $^L G$ \cite{L}. 
This parabolic subgroup corresponds to a maximal parabolic subgroup
 of $SO(2m+1,1)$
 whose Levi subgroup $L$ is a real form of $SO(2m,\bC) \times SO(2,\bC)$ and thus is isomorphic to $ SO(2m,0) \times SO(1,1)$. 
Note that $SO(1,1)$ is a disconnected group.  

By \cite[p.~35]{GGP} there are $2^m$ representations in the Vogan packet containing $\Pi_{m,\delta}$
 and they are parametrized by the characters
 of a finite group ${\mathcal{A}}_1 \simeq (\bZ_2)^m$.  
We write $VP(\Pi_{m,\delta})$
for this Vogan packet.

We can describe the representations  in the  Vogan packet $VP(\Pi_{m,\delta})$ as  follows: 
for $l+k=2m+2$, 
 a real form $SO(l,k)$ is called {\it{pure}} if $l$ is odd. 
The Levi subgroups of   parabolic subgroups  in the same Vogan packet are isomorphic to $L=SO(2m-2p,2p)\times SO(1,1).$ Principal series representations, which are induced from the outer tensor product of a discrete series representation with the same infinitesimal character
 as $\Exterior^{\frac n 2 }({\mathbb{C}}^{n})_+$
 (or $\Exterior^{\frac n 2 }({\mathbb{C}}^{n})_-$)
 and a one-dimensional representation of $SO(1,1)$, 
 are irreducible. 
These induced representations  are in  $VP(\Pi_{m,\delta} )$ 
 if they have the same central character as $\Pi_{m,\delta}$  \cite{V}.

\medskip
\noindent
\begin{remark}
\begin{itemize}
\item[{\rm{(1)}}]
The Vogan packet containing $(SO(1+2m,1),\Pi_{m,\delta})$, $\delta= (-1)^{m+1}$ does contain neither the pair $(SO(2m+2,0)$, finite-dimensional representation) nor ($SO(0,2m+2)$, finite-dimensional representation).  
\item [{\rm{(2)}}]
If $(SO(1+2m,1),\Pi)$ is in $VP(\Pi_{m,\delta})$,
 then $\Pi = \Pi_{m,\delta}$.  

\item[{\rm{(3)}}] 
By abuse of notation
 we say that the Vogan packet containing $(SO(1+2m,1),\Pi_{m,\delta})$
 also contains $(SO(2m+1,1),\Pi_{m,\delta})$.  

\item[{\rm{(4)}}]
Using the same considerations for $m-1$ and $\delta = (-1)^m$
 we obtain a Vogan packet $VP(\varpi_{m-1,\delta} )$
 which contains the pair
 $(SO(2m-1,1), \varpi_{m-1,(-1)^m})$. 
\end{itemize}
\end{remark}

\subsection{Gross--Prasad conjecture I:
Symmetry breaking from  $\Pi_{m,-}$  to the discrete series representation $\pi_m$}
\label{subsec:IV.4}

We consider the Vogan packet of tempered representations of $SO(2m+1,1)\times SO(2m,1)$ which contains the pair $(SO(2m+1,1) \times SO(2m,1), \Pi_{m,\delta} \boxtimes \pi_m )$, 
 or the Vogan packet which contains the pair
$(SO(1,1+2m)\times SO(1,2m),  \Pi_{m,\delta} \boxtimes \pi_m )$. 
The representations in these packets are parametrized by characters of
\[
     {\mathcal{A}}_1 \times {\mathcal{A}}_2 
     \simeq 
     (\bZ_2)^{\color{black}{m}} \times (\bZ_2)^m \simeq (\bZ_2)^{2m} . 
\] 
B.~Gross and D.~Prasad propose an algorithm
 which determines a pair 
$
     \chi_1\in \widehat{{\mathcal{A}}_1},\chi_2 \in \widehat{{\mathcal{A}}_2}$
 hence representations 
\[
   (\Pi(\chi_1),\pi(\chi_2 )) \in VP(\Pi_{m,\delta}) \times VP(\pi_m)
\]
 so that 
\[ 
   \mbox{Hom}_{G(\chi_2)}(\Pi(\chi_1)|_{G(\chi_2)},\pi(\chi_2) ) \not = \{0\},
\]
where $G(\chi_2) $ is the pure inner form determined by $\chi_2$.

\medskip
Let $T_\bC$ be a torus in $SO(2m+2,\bC)  \times SO(2m+1,\bC)$, 
 and $X^{\ast}(T_\bC)$ the character group. 
Fix basis 
\[
   X^*(T_{\mathbb{C}})
   = 
   \bZ e_1\oplus \bZ e_2 \oplus \dots \oplus \bZ {{e_{m+1}}} \oplus \bZ f_1\oplus \bZ f_2 \oplus \dots \oplus \bZ f_m  
\]
such that the standard root basis $\Delta_0$ is given by
\[
e_1-e_2, e_2-e_3,\dots,{{e_{m}-e_{m+1},e_{m}+e_{m+1}}}, f_1-f_2 ,f_2-f_3, \dots, f_{m-1}-f_m,f_m
\]
{if $m \ge 1$}.

We   fix as before  $\delta=(-1)^{m+1}$.  
We can identify the Langlands parameter of the Vogan packet containing 
\[
(SO(2m+1,1)\times SO(2m,1), \Pi_{m, {{\delta}}} \boxtimes \pi_m)
\]
with 
 \[ me_1+(m-1)e_2+ \dots +e_m+0e_{m+1} +(m-\frac 1 2 )f_1+(m- \frac 3 2 ) f_2 +\dots + \frac 1 2 f_m.\]

Let $\delta_i$ be the element which is $-1$ in the $i$th factor
 of ${\mathcal{A}}_1$
 and equal to $1$ everywhere else
 and $\varepsilon_j$ the element which is $-1$
 in the $j$th factor of ${\mathcal{A}}_2$ and $1$ everywhere else.

Then the algorithm \cite[p.~993]{GP} determines
 $\chi_1 \in \widehat {{\mathcal{A}}_1}$
 and $\chi_2 \in \widehat {{\mathcal{A}}_2}$
 by 
\[
 \chi_1(\delta_i)= (-1)^{\#m-i+1>}  \ \ \ \ 
 \mbox{and } 
 \chi_2(\varepsilon_j)=  (-1)^{\#m-j + \frac 1 2<}, 
\]
where $\#m-i+1>$ is the cardinality
 of the set $\{j: m-i+1> \text{the coefficients of $f_j$}\}$,
 and $\# m-j+\frac 1 2 <$ is the cardinality
 of the set $\{i:m-j+\frac 1 2 < \text{the coefficients of $e_i$}\}$.

We normalize the quasi-split form by 
\begin{alignat*}{2}
 G^o &= SO(m+1,m+1) \times SO(m,m+1) \quad &&\mbox{if $m$ is even,} 
\\
 G^o &= SO(m+2,m) \times SO(m+1, m)  \quad &&\mbox{if $m$ is odd.} 
\end{alignat*}
Applying the formul{\ae} in \cite[(12.21)]{GP}
 we define the integers $p$ and $q$
 with $0 \leq p \leq m $ and $0 \leq q \leq m$
by
\[
 p  =  \# \{ i : \chi_1(    \delta_i  )  = (-1)^i  \}  \quad \mbox{ and }  
\quad  q  =  \# \{ j:  \chi_2(   \varepsilon_j     )  = (-1)^{m+j}  \}
\]
and we get the pure forms
\[G = SO(2m-2p+1, 2p+1) \times SO(2q ,2m-2q+1) \quad \mbox{if $m$ is even,}\]
\[G= SO(2p+1,2m-2p+1) \times SO(2m-2q, 2q+1) \quad \mbox{if $m$ is odd.}\]
In our setting,  we get the pair of integers $(p,q)=(0,m)$ for $m$ even; $(p,q)=(m,0)$ for $m$ odd.
Applying \cite[(12.22)]{GP} with correction by changing $n$ by $m$ 
 {\it{loc.cit.}},
we  deduce that this character defines the pure inner form 
\begin{itemize}
\item $G=SO(2m+1,1) \times SO(2m,1)$ if $m$ is even,
\item $G=SO(2m+1,1) \times SO(2m, 1)$ if $m$ is odd.
\end{itemize}
The only representation in $VP(\Pi_{m,\delta} )\times VP(\pi_m)$ with this pair of pure inner forms  is $\Pi_{m,\delta} \times \pi_m$.
Hence Theorem \ref{symbreaking/rho} implies the following.  

\medskip

\noindent
{\bf  Conclusion:} {\it The result 
\[
   \dim_{\mathbb{C}}
   {\operatorname{Hom}}_{SO(2m,1)}(\Pi_{m,{(-1)^{m+1}}}|_{SO(2m,1)},\pi_m)  =1
 \] 
confirms the conjectures by B.~Gross and D.~Prasad } \cite{GP}.

\medskip

\subsection{Gross--Prasad conjecture II: Symmetry breaking
 from the discrete series representation $\pi_m$ to ${\varpi}_{m-1,{(-1)^{m}}}$}
\label{subsec:GPconjII}

We consider  the Vogan packet of tempered representations containing the pair $(SO(2m,1) \times SO(2m-1,1), \pi_{m} \boxtimes \varpi_{m-1,(-1)^m} )$, 
 {\it{i.e.}}, 
 the Vogan packet 
\[
VP(\pi_m\boxtimes \varpi_{m-1,(-1)^m})\subset VP(\pi_m) \times VP(\varpi_{m-1,(-1)^m}).  
\] 
The packet $VP(\pi_m) \times VP(\varpi_{m,(-1)^m})$ is
 parametrized by  characters of the finite group  
\[
     {\mathcal{A}}_2 \times {\mathcal{A}}_3
 \simeq 
 (\bZ_2)^m \times (\bZ_2)^{m-1}\simeq (\bZ_2)^{2m-1}.
\]
As in Section \ref{subsec:IV.4} 
 we use the algorithm by B.~Gross and D.~Prasad to determine a pair 
 $\chi_2 \in \widehat{{\mathcal{A}}_2}, \chi_3 \in \widehat{{\mathcal{A}}_3}$
 and hence representations 
\[
 (\pi(\chi_2), \varpi(\chi_3 )) \in VP(\pi_{m})\times VP(\varpi_{m-1,(-1)^m}  )\] so that 
\[ 
     \mbox{Hom}_{G(\chi_3)}(\pi(\chi_2)|_{G(\chi_3)},\varpi(\chi_3) ) \not = \{0\}, 
\]
where $G(\chi_3) $ is the pure inner form determined by $\chi_3$.

\medskip

Let $T_\bC$ be a torus in $SO(2m+1,\bC)\times  SO(2m, \bC)$
 and $X^{\ast}(T_{\mathbb{C}})$ the character group.  
Fix basis 
\[
X^*(T_\bC)
= \bZ f_1\oplus \bZ f_2 \oplus \dots 
\oplus \bZ f_m \oplus \bZ g_1\oplus \bZ g_2 \oplus \dots \oplus \bZ {g_{m}} 
\]
such that the standard root basis $\Delta_0$ is given by
\[
f_1-f_2, f_2-f_3,\dots,f_{m-1}-f_m,f_m, g_1-g_2 ,g_2-g_3, \dots, 
 {g_{m-1}-g_{m},g_{m-1}+g_{m}}
\]
{for $m \ge 2$}.

Fix as before  $\epsilon=(-1)^{m}$.  
We identify the Langlands parameter of the Vogan packet 
\[
VP( \pi_{m} ) \times VP( \varpi_{m,\epsilon})
\]
 with 
 \[ (m-\frac 1 2 )f_1+(m- \frac 3 2 ) f_2 +\dots + \frac 1 2 f_m +(m-1)g_1+(m-2)g_2+ \dots +g_{m-1}+0g_{m} .\]

Again applying \cite[Prop.~12.18]{GP}
 we define characters $\chi_2 \in \widehat{{\mathcal{A}}_2} $,
 $\chi_3 \in \widehat{{\mathcal{A}}_3}$ as follows: 
Let $\varepsilon_j \in {\mathcal{A}}_2 \simeq (\bZ_2)^m$ be
 the element which is $-1$ in the $j$th factor
 and equal to $1$ everywhere else
 as in Section \ref{subsec:IV.4};
 $\gamma_k \in {\mathcal{A}}_3 \simeq (\bZ_2)^{m-1}$ the element which is $-1$ in the $k$th factor and $1$ everywhere else.
Then $\chi_2 \in \widehat{{\mathcal{A}}_2}$ and 
 $\chi_3 \in \widehat{{\mathcal{A}}_3}$ are determined by 
\[ \chi_2(\varepsilon_j)= (-1)^{\#m- j +1/2<}  \ \ \ \ \mbox{and }
 \chi_3(\gamma_k)=  (-1)^{\#m-k >}\]
where $\#m-j+\frac 1 2 <$ is the cardinality
 of the set $\{k:m-j+\frac 1 2 < \text{the coefficients of $g_k$}\}$,
 and $\# m-k >$ is the cardinality of the set
 $\{j:m-k > \text{the coefficients of $f_j$}\}$.

We normalize the quasi-split form by 
\begin{alignat*}{2}
 G^o &= SO(m+1, m) \times SO(m+1,m-1) \quad &&\mbox{if $m$ is even,} 
\\
 G^o&= SO(m,m+1) \times SO(m , m)     \quad &&\mbox{if $m$ is odd.} 
\end{alignat*}

Applying the formul{\ae} in \cite[(12.21)]{GP}
 we define the integers $p$ and $q$
 with $0 \leq p \leq m$ and $0 \leq q \leq m-1$
by
\[
 p  =  \# \{ j : \chi_2(    \varepsilon_j  )  = (-1)^j  \}  
\quad \mbox{ and }  
\quad  q  =  \# \{ k :  \chi_3(   \gamma _k      )  = (-1)^{m+k}  \}, 
\]
and we get
\[G = SO(2m-2p+1, 2p) \times SO(2q+1 ,2m-2q-1) \quad \mbox{ if $m$ is even,}\]
\[G= SO(2p+1,2m-2p) \times SO(2m-2q-1, 2q+1) \quad \mbox{ if $m$ is odd.}\]

In our setting, the pair of integers 
 $(p,q)$ is given by $(p,q)=(m,0)$ for $m$ even;
 $(p,q)=(0,m-1)$ for $m$ odd.
Applying \cite[(12.22)]{GP} 
we  deduce that this character defines the pure inner form 
\begin{itemize}
\item $G^0=SO(1,2m) \times SO(1,2m-1)$ if $m$ is even,
\item $G^0=SO(1,2m) \times SO(1, 2m-1)$ if $m$ is odd.  
\end{itemize}

The only representation in $VP(\pi_{m}) \times VP(\varpi_{m-1,\epsilon})$ with this pair of pure inner forms  is $\pi_{m} \times \varpi_{m-1,\epsilon}$.

\medskip

\noindent
{\bf  Conclusion:} {\it The result 
\[
   \dim_{\mathbb{C}}
   {\operatorname{Hom}}_{SO(1,2m-1)}(\pi_{m}|_{SO(1,2m-1)},\varpi_{m-1,(-1)^{m}})  =1
 \] 
confirms the conjectures by B.~Gross and D.~Prasad } \cite{GP}.

\medskip
{\bf{{Acknowledgements}}}\enspace
The first author was partially supported by Grant-in-Aid for Scientific 
Research (A)
(25247006), Japan Society for the Promotion of Science.  

Research by B.~Speh was partially supported by  NSF grant DMS-1500644. 
Part of this research was conducted during a visit of the second author at  the Graduate School of  of Mathematics of the University of Tokyo, Komaba. 
She would like to thank  it for its support and hospitality during her stay.

\bigskip

\color{black}


\begin{thebibliography}{M}

\bibitem{ABV} J.~Adams, D.~Barbasch, D.~Vogan, 
The Langlands Classification and Irreducible Characters for Real Reductive Groups, 
 Progr. Math.  {\bf{104}}, 
Birkh\"{a}user, Boston--Basel--Berlin, 1992.
 

\bibitem{BW} A.~Borel and N.~Wallach, 
Continuous Cohomology, Discrete Subgroups, and Representations of
Reductive Groups. 2nd edition. Mathematical Surveys and Monographs, 
{\bf{67}}, Amer.~ Math.~ Soc., Providence, 2000.  

\bibitem{GGP} W.-T.~Gan, B.~Gross, D.~Prasad: Symplectic local root numbers, central critical $L$-values, and restriction problems in the representation theory of classical groups, Ast\'{e}risque {\bf{346}}, 1--109, 2012.


\bibitem{GR} B.~Gross and M.~Reeder, 
{From Laplace to Langlands via representations of orthogonal groups},
Bull. Amer. Math. Soc. (N.S.), 
{\bf{43}},  
(2006), 
163--205.  



\bibitem{GP} B.~Gross and D.~Prasad, 
On the decomposition of a representation of $\operatorname{SO}_n$
 when restricted to $\operatorname{SO}_{n-1}$, 
Canad. J. Math. {\bf{44}}, 
 (1992), no. 5, 
974--1002.  

  

\bibitem{GW}  
B.~Gross and N.~Wallach, 
 Restriction of small discrete series representations to symmetric subgroups, Proc. Sympos. Pure Math. {\bf{68}},
 (2000), Amer.~Math.~Soc., 255--272.  

\bibitem{KV} A.~W.~ Knapp and D.~ Vogan, 
Cohomological induction and unitary representations.
Princeton Mathematical Series, {\textbf{45}}. 
Princeton University Press, Princeton, NJ, 1995. xx+948 pp. 
ISBN: 0-691-03756-6.  

\bibitem{KMemoirs92}
T.~Kobayashi, 
Singular Unitary Representations and Discrete Series
 for Indefinite Stiefel Manifolds $U(p,q;F)/U(p-m,q;F)$, 
Mem.~Amer.~Math.~Soc.~ 
\href{http://www.ams.org/books/memo/0462/}{{\bf{462}}}, Amer. Math. Soc., 1992. 106 pp. 
ISBN:9810210906.  

\bibitem{KInvent98}
T.~Kobayashi, 
Discrete decomposability of the restriction of $A_{\mathfrak{q}}(\lambda)$
 with respect to reductive
subgroups. III. Restriction of Harish-Chandra modules and associated varieties, Invent. Math.
{\bf{131}}, (1998), no. 2, 
\href{http://dx.doi.org/10.1007/s002220050203}
{229--256}.  

\bibitem{KKP} T.~Kobayashi, T.~Kubo, M.~Pevzner, 
Conformal Symmetry Breaking Differential Operators
 for Differential Forms on Spheres, 
Lecture Notes in Math., 
 \href{http://dx.doi.org/10.1007/978-981-10-2657-7}{vol.~{\bf{2170}}}, 
ix $+$ 192 pages, 
 2016. 
ISBN: 978-981-10-2657-7.  

\bibitem{KO} 
T.~Kobayashi and T.~Oshima, 
Finite multiplicity theorems for induction and restriction, Adv.
Math. {\bf{248}}, (2013), 
\href{http://dx.doi.org/10.1016/j.aim.2013.07.015}
{921--944}. 

\bibitem{KS} T.~Kobayashi and B.~Speh, 
Symmetry breaking for representations of rank one orthogonal groups, 
Mem. Amer. Math. Soc., 
vol. 
\href
{DOI: http://dx.doi.org/10.1090/memo/1126}
{\bf{238}}, 
No. 1126, 
(2015), 
v $+$ 112 pages, 
ISBN: 978-1-4704-1922-6.  

\bibitem{L} R.~Langlands, On the classification of irreducible representations of real reductive groups, Math. Surveys and Monographs 
{\bf{31}}, 
Amer.~Math.~Soc., 
Providence,
1988.  

\bibitem{SunZhu}
B.~Sun and C.-B.~Zhu, 
Multiplicity one theorems: the Archimedean case, Ann. of Math. (2), 
{\bf{175}}, (2012), no. 1, 
\href{DOI 10.4007/annals.2012.175.1.2}
{23--44}. 

\bibitem{V}
D.~Vogan, 
The local Langlands conjecture, 
Contemp. Math.,  {\bf{145}}, (1993), 305--379,
 Amer. Math. Soc. 

\bibitem{VZ} 
D.~Vogan and G.~Zuckerman, 
Unitary representations with nonzero cohomology, 
Compositio Math. 
{\bf{53}}, 
(1984), 
no. 1, 51--90.  

\bibitem{W} 
N.~R.~Wallach, 
Real Reductive Groups. II, 
Pure and Applied Mathematics, 
{\bf{132}}, Academic Press, Inc., Boston, MA, 1992. 
ISBN 978-0127329611.  

\end{thebibliography}
\end{document}